\documentclass[a4,11pt,authoryear]{elsarticle}
\usepackage[margin=1in]{geometry}
\usepackage[utf8]{inputenc}

\usepackage{float} 
\usepackage{soul}
\usepackage{multirow}
\usepackage{graphicx}
\usepackage{outlines}
\usepackage{longtable}
\usepackage{tabu}
\usepackage{footnote}
\usepackage{mathtools}
\usepackage{amssymb} 
\usepackage{blkarray, bigstrut}
\usepackage[linesnumbered,lined,boxed]{algorithm2e}
\usepackage{setspace}
\usepackage{comment}
\usepackage{url}

\usepackage{breakurl}

\usepackage{natbib}
 \bibpunct[, ]{(}{)}{,}{a}{}{,}%
\makeatletter
\def\ps@pprintTitle{%
 \let\@oddhead\@empty
 \let\@evenhead\@empty
 \def\@oddfoot{\centerline{\thepage}}%
 \let\@evenfoot\@oddfoot}
\makeatother

\journal{European Journal of Operational Research}

\begin{document}

\begin{frontmatter}


\title{Decision Programming for optimizing the Finnish Colorectal Cancer Screening Program} 

\author[BIZ]{Lauri Neuvonen\corref{cor1}}\ead{lauri.neuvonen@aalto.fi}
 \cortext[cor1]{Corresponding author.}
\author[BIZ]{Mary Dillon}
\author[BIZ]{Eeva Vilkkumaa}
\author[SCI]{Ahti Salo}
\author[FCR]{Maija Jäntti}
\author[FCR]{Sirpa Heinävaara}
\address[BIZ]{Department of Information and Service Management, Aalto University School of Business, Finland}
\address[SCI]{Department of Mathematics and Systems Analysis, Aalto University School of Science, Finland}
\address[FCR]{The Finnish Cancer Registry}

\begin{abstract}
  In Finland colorectal cancer (CRC) incidence rates have steadily increased over the last decades and as of 2017, CRC is the sixth most common cause of death. CRC is a crucial concern for the public health of Finland. We optimize the faecal immunochemical test (FIT) cut-off level for specified target populations with regard to minimization of cancer prevalence with Decision Programming, which is a novel approach to solving discrete multi-stage decision problems under uncertainty. The results present optimal cut-off levels for Finnish target groups with different colonoscopy capacity constraints. 
  Finally, the estimated resulting cancer prevalence, the amount of required colonoscopies and third-party payer costs resulting from found strategies to the newly implemented Finnish CRC Screening Programme, which have not previously been published, are presented.  
\end{abstract}


\end{frontmatter}



\section{Introduction}


Colorectal cancer (CRC), also referred to as bowel cancer, is cancer of the colon and/or rectum. CRC incidence rates have steadily increased over the last decades with CRC currently being the second most common cancer in adults in Finland, and the third most common worldwide \citep{FCR_cancer,WHO_2018}. As of 2017, CRC is the sixth most common cause of death, whereby it is clearly a crucial concern to public health.

 Over 70\% of colorectal cancers develop via the adenoma-carcinoma sequence \citep{hardy2000molecular} in which adenomas (i.e., growths  on the epithelial tissue) develop into cancer. This sequence is a slow process that can take from several years to a decade \citep{simon2016colorectal}. Due to this, early detection and removal of pre-cancerous adenomas can prevent progression to cancer, making CRC extremely suitable for population screening. Population screening of CRC using faecal occult blood testing has long been shown to reduce later stage CRC incidence and improve mortality rates \citep{mandel1993reducing,kewenter1994results,kronborg1996randomised,mandel1999colorectal}. Along with early detection of adenomas, early diagnosis of cancer improves the probability of successful treatments and outcomes by providing care at the earliest possible time. This has an important impact on public health strategy as it helps avoid CRC-attributed deaths and disabilities, in addition to high treatment and indirect costs associated with advanced cancer stages \citep{WHO_2017}.
 
 In April 2019, Finland established a new CRC screening programme in volunteering municipalities. The programme will become nationwide in 2022, and by 2031 all people aged between 56 and 74  will be invited to screen. In the new programme, the feacal immunochemical test (FIT) replaces the previously used guaiac-based faecal occult blood test (gFOBT), as the European guidelines for quality assurance in CRC screening recommend FITs over gFOBTs due to higher sensitivity levels and increased participation rates \citep{halloran2012european,allison2014population, bretagne2019switching}. The programme invites entire age cohorts to screen using a FIT, after which those participants with a positive FIT result are invited to a colonoscopy. 
 
 In the current programme, the FIT cut-off level (hemoglobin level in stool sample) for a positive result is fixed based on the participant's sex.  In particular, the Finnish programme uses a cut-off level of 25 mg or ng/g for females and 70 mg or ng/g for males, regardless of their age. In reality, however, the interpretation of the result should be different for participants of different risk profiles (i.e., age, sex, and family history of cancer). Moreover, because a positive result leads to an invitation to a colonoscopy, the cut-off level should also take into account the capacity of the participating municipality to carry out colonoscopies. Consequently, the programme may not be the most cost-effective approach to reduce Finnish colorectal cancer incidence and mortality. 

Traditionally, cost-effectiveness evaluations of public health programmes, such as population screening, are assessed through methods of cost-benefit analysis (e.g.~\citealp{ellison2002cost}), cost-utility analysis (e.g.~\citealp{gupta2011endoscopy}), or cost-effectiveness analysis (e.g.~\citealp{dillon2018family}). Health economic methodologies like these typically compare costs and health outcomes to a baseline strategy (for example, no screening). Quality of life indicators, such as quality-adjusted life-years (QALYs) and life-years saved/gained (LYS/LYG), are frequently used to calculate the incremental cost-effectiveness ratio (ICER) for comparing alternative strategies. In these analyses, strategies are not optimized in the mathematical sense, but rather assessed in regard to their dominance. Consequently, such analyses tend to suggest strategies which are either unaffordable due to a lack of budget constraint or suboptimal in that (i) the resources could be reallocated to achieve a better population-level health outcome or (ii) the same health outcome could be attained with less resources.

To identify optimal strategies within cost-effectiveness analyses, decision analytic models such as decision trees have been employed (e.g.,~\citealp{hynninen2019value},~\citealp{hynninen2021operationalization}). These models are typically static in that they do not account for possible changes in a participant's state of health over time, whereby they cannot be readily applied to support the optimization of screening programmes. To capture time-dependence in a participant’s state of health, several models utilizing Markov Decision Processes (MDPs) have been proposed. In MDP models, the participant is assumed to move from one state to another (e.g., from being healthy to having an adenoma or cancer) with state transition probabilities that can be affected by carrying out preventive or treatment actions. In Partially Observable Markov Decision Process (POMDP) models in particular, the state of a participant may not be perfectly observed but information regarding this state can be obtained through screening, for instance. \citet{Ayer2012} build a POMDP model to optimize participant-specific mammography screening times. \citet{alagoz2013optimal} provides a tutorial in optimizing cancer screening using a POMDP approach and \citet{Erenay2014OptimizingSurveillance} develop a finite-horizon POMDP model to determine optimal CRC screening policies.  Whereas, \citet{cevik2018analysis} present a constrained POMDP model to study the optimal allocation of limited mammography resources to screen a population. \citet{lee2019optimal} optimize the use of limited resources for the screening of a population for hepatocellular carcinoma by modeling the problem as a family of restless bandits in which each participant’s disease progression is assumed to evolve as a POMDP.

POMDP approaches are very useful when the possible deterioration in a participant’s state of health over time needs to be accounted for. Nevertheless, POMDP studies often assume fixed interpretations of test results (i.e., predetermined positivity cut-off levels). In particular, few of the above studies recognize that the optimal interpretation of test results may depend on how much resources there are. Moreover, the above studies employ dynamic programming to solve the optimization problems. For this reason, the models can only have a single objective such as maximizing the expected health outcome or expected net benefit, wherein health outcomes are converted into monetary terms using some parameter such as the willingness to pay threshold (e.g.,~\citealp{chen2018optimal}). Importantly, focusing on expected outcomes fails to accommodate risk considerations, which could be highly relevant in the context of cancer screening. 

In this paper, we optimize population screening of CRC in regard to minimizing cancer prevalence, i.e. share of population with CRC, in the target population with regard to the maximum number of colonoscopies that can be performed in a year. Toward this end, we build a two-level optimization model. At the first level, we find all Pareto optimal screening strategies for each segment of the population, defined by the participants' age and sex, maximizing the probabilities of discovering a cancer or a large growth and minimizing the probability of performing a colonoscopy. At the second level, we identify the combination of segment-specific strategies that result in the highest expected number of detected cancers, given a colonoscopy constraint. 
In building this model, we employ Decision Programming \citep{salo2022decision}, a novel approach to solving discrete multi-stage decision problems under uncertainty. In contrast to earlier approaches to optimizing screening programmes, Decision Programming (i) finds all Pareto optimal solutions with respect to multiple objectives, and (ii) permits the use of deterministic and probabilistic risk measures as objectives or constraints in the model formulation. As a methodology it has not been applied to the health care context. 

The contribution of this paper to existing literature is twofold. First, we propose defensible decision recommendations for developing the CRC screening program in Finland. Second, we demonstrate the usefulness of the Decision Programming methodology in health care context.

The paper is structured as follows: Section 2,  Methods, including case study and model description; Section 3 details the decision programming model formulation; Section 4 presents the  case study results; Section 5 discusses the model, the results and the conclusions including strengths and limitations of the method and future research directions.

\section{Problem overview}\label{Sec:prob_overview}

The new colorectal cancer screening programme in Finland commenced in April 2019 (see~\citealp{CanReg2019} for details). In 2021 all 60-, 62-, 64-, 66- and 68-year-old males and females in volunteering municipalities are invited to participate in the programme. Presently, twelve out of 311 municipalities in Finland have volunteered to take part in the programme. The programme will become nationwide in 2022.

A central screening hub mails the screening invitations and FITs, analyzes the samples, and mails the result letters to participants. The FIT is performed at home and returned in a pre-paid envelope to the screening hub. Participation in the screening programme is free of charge for the invitee. Once the laboratory has analyzed the FIT sample, a written result is given to the participant via mail. Those with a positive result are required to contact a screening nurse in their municipality of residence to discuss the need for further examination. This additional examination is usually a colonoscopy; however, this may vary between municipalities depending on their standard procedure. If a growth is found during the additional examination, a sample is taken and sent to a pathologist for analysis. Once the pathologist's findings have been reviewed, a decision on the need for further examination or treatment is made. If necessary, the participant is referred to surgery. Treatment may be required at this stage. However, the treatment process is not within the scope of the screening programme. After the screening period has ended, i.e. two-years have passed since the invitation, a new invitation is sent and the next screening period commences. Screening is continued in this periodic manner as long as the invitee's age is within the program limits.

The main decisions affecting the outcomes of the screening program include what population segments are invited, what pre-screening and examination methods are used. Depending on these selections, a certain number of examinations will be performed, costs generated and health benefits gained. However, there's only a limited amount of colonoscopies that can be performed in a year. This gives rise to the task of deriving maximum health benefits from these resources, i.e. maximizing screening efficiency. 

This paper develops a model to produce optimal screening strategies for the Finnish colorectal cancer screening programme. The overall objective is to minimize the prevalence of cancers while minimizing the total colonoscopies performed. A screening strategy is defined as a combined selection of decisions for all target population segments. For this problem, we include the FIT cut-off level, whether to offer an incentive, and whether to send an invite 
as decisions. In the Finnish programme, the FIT cut-off levels are fixed at 70 \textit{$\mu$g/g} or 410 \textit{ng/ml} for males and 25 \textit{$\mu$g/g} or 150 \textit{ng/ml} for females, (see \citet{shaukat2015colorectal} Table 4.3 for unit conversion). Please note that the cutoff level stated by the manufacturer is 100 \textit{ng/ml} \citep{shaukat2015colorectal}, however the FIT can measure up to 1000 \textit{ng/ml} \citep{FOBgold}. The inclusion of an incentive decision in our problem does not reflect the reality of the programme; however, it may be an option in other countries with similar programmes.

\section{Optimisation of the screening programme.}

The overall objective is to find a strategy, i.e. a combination of invitation, FIT cut-off level and incentive decisions for each segment, that minimises cancer prevalence in the target population with a limited colonoscopy resource. We model the problem of finding such a strategy via two phases. We first present an overview of the phases to create context, and later the details in dedicated sections. 

\begin{figure}[!ht]
    \centering
    \includegraphics[width=\textwidth]{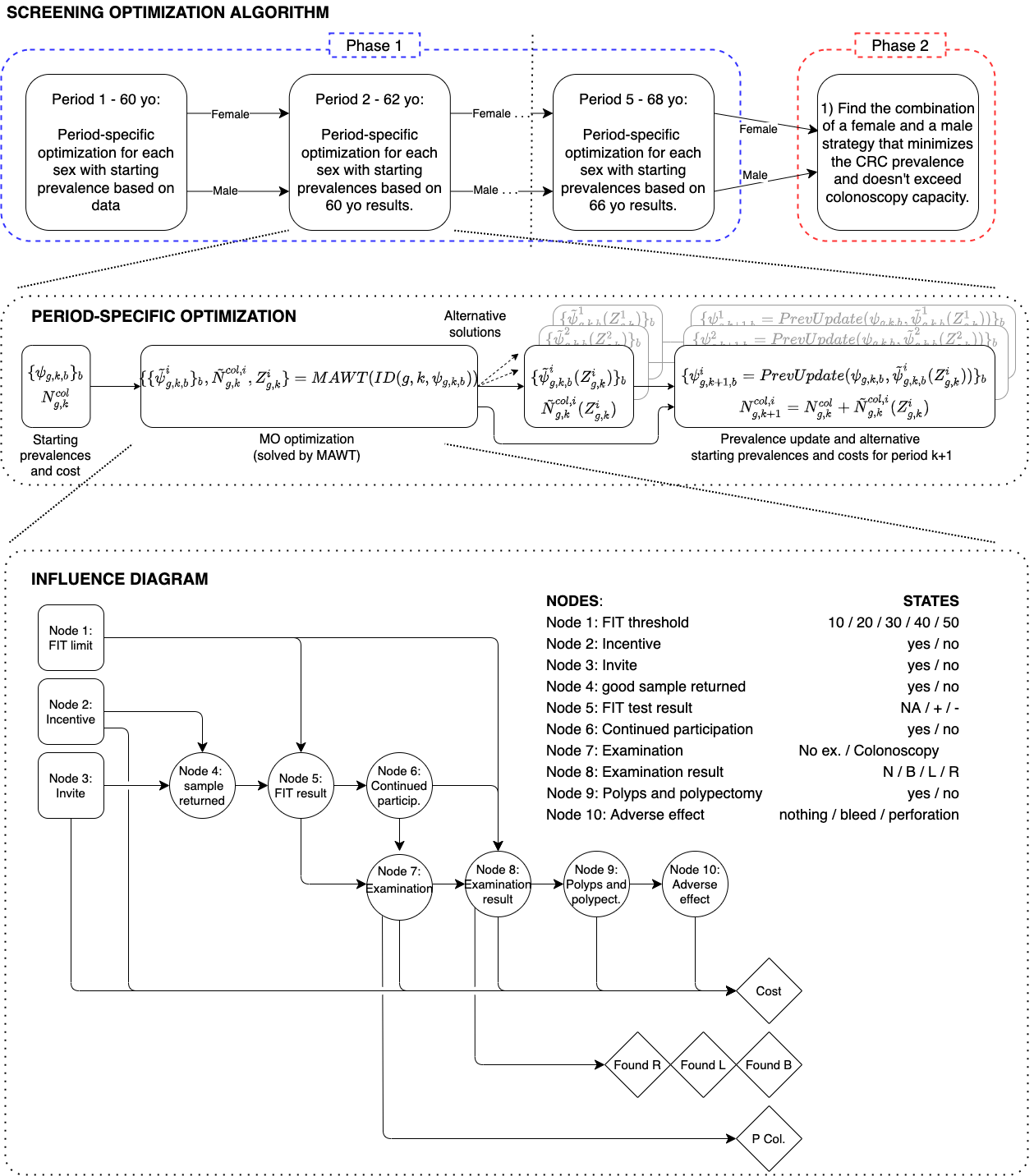}
    \caption{Schematic description of the multistage optimization approach to solve the CRC screening optimization problem.}
    \label{fig:combined_algo}
\end{figure}

In the first phase, we generate the set of Pareto optimal strategies $Z^*_{g}$ for each sex ($g$) for all periods ($k$), i.e. segments of the population, that minimizes the expected number of colonoscopies, while maximizing the expected number of detected growths and cancers in each screening round. In the second phase, we combine the sex-specific strategies to create an overall strategy $Z$ that minimizes the prevalence of cancer in the screened population within a given resource constraint. Figure \ref{fig:combined_algo} illustrates our two phased approach. 

The first phase begins with identifying the set of optimal screening strategies for those in period $k=1$, i.e. aged 60, of separate sex. Each period represents one screening round, therefore, in the following period (i.e. two years later) optimal strategies for those aged 62 are added to the set. This continues until the final period is reached and we have a set of optimal screening strategies for both sexes and all ages included in the screening programme. As the solution at $k$ is dependent on the prevalence distribution produced by the preceding solution at $k-1$, the segment specific solutions form a tree structure, each Pareto optimal solution creating a new branch. We call a path through this tree, i.e. a collection of segment-specific strategies, one for each period, a \emph{sex-specific strategy} $Z_{g}$. From here, the second phase identifies a sex-specific strategy for all-aged males and for all-aged females from the final set of sex-specific efficient strategies. This latter phase is a simple selection problem. However, Phase 1 is a relatively complex optimization problem.

\subsection{Main assumptions}


We assume that the sex-specific strategies are stable in the sense that the same strategy is applied to e.g. 62-year old males regardless of the year the entered the screening programme. This means that the strategy applied to one group of people during their whole screening programme is the same that is applied to the different age groups within one screening period. This reasonable assumption simplifies the problem of finding sex-specific strategies significantly. 
The population segments are assumed to be internally homogenous, i.e. there's no difference between what is known about two 62-year old males before they enter the program. 

The transitions from different bowel to other states between beginnings of periods are assumed to be linear. This means that during a period a portion of people with normal bowels will develop benign growths, a portion of people with benign growths will develop large adenomas etc. In addition, we assume that all participants for whom a abnormal bowel state is identified, are treated and their bowel state becomes normal. These rather strong assumptions have been discussed with the Finnish cancer registry and can be deemed acceptable for the purposes of this paper. However, more refined transition models could be integrated into this modeling approach.


\subsection{Segment specific model}
\label{sec:ID_MAWT}

The segment specific efficient strategies are generated by combining the \emph{Modified Augmented Weighted Tchebychev} (MAWT) norm approach \citep{Holzmann2018SolvingScalarizations} with an Influence Diagram (ID) capturing the screening process for one participant in a segment. The influence diagram has been modeled in Decision Programming style \citep{salo2022decision} to allow inclusion of risk measures. Combining these elements leads to the optimisation problem described in Equations \eqref{eq:MAWT_obj}-\eqref{eq:pi_def} that can be solved using the MAWT algorithm described in \cite{Holzmann2018SolvingScalarizations}. The differences between target segments are reflected in the parameter values of the IDs. In particular, the probability of a participant's bowel to be in a certain state $b$ is the prevalence $\psi_{g,k,b}$ of that bowel state in the target segment $(g,k)$ the participant belongs to.

The Finnish screening programme described in the previous section can be captured by an influence diagram, as shown in Figure \ref{fig:combined_algo}. An influence diagram is a discrete acyclic graph constructed of three types of nodes $\mathcal{N}$ represented by sets; decision $\mathcal{D}$, chance $\mathcal{C}$ and values nodes $\mathcal{V}$, with dependencies shown by directed arcs $\mathcal{A} \subseteq \{(i,j)|i,j \in \mathcal{N}, i \neq j\}$. In Figure \ref{fig:combined_algo}, decisions in the screening pathway are represented by squares, uncertainties (or chances) in the screening process are represented by circles, and the values that are to be optimized are represented by diamonds. Every chance $j \in \mathcal{C}$ and decision $j \in \mathcal{D}$ node has a finite set of discrete states $s_j \in \mathcal{S}_j$. The state $s_j \in \mathcal{S}_j$ of a given node $j \in \mathcal{C} \cup \mathcal{D}$ represents a chance or decision alternative, a full list of these alternatives for our problem are presented below. An arc $(i,j) \in \mathcal{A}$, represented by an arrow, indicates that node $i$ is the predecessor of node $j$, and that the state at node $j$, $s_j$, is conditionally dependent on the state at the preceding node, $s_i$. 


In Figure \ref{fig:combined_algo}, the decision nodes 1, 2, and 3 correspond to decisions about 1) what FIT test cut-off value to use to select participants for a colonoscopy in the given target segment, 2) whether to use an incentive to boost participation rate among invitees in this segment (specifically, we   assume that an incentive worth 50 euros halves the number of non-returned samples) and 3) whether to invite the target segment to the screening program. Let us denote by $s_j\in \mathcal{S}_j, j\in \mathcal{D}$ the alternatives for decision $j$. Here, the sets  $\mathcal{S}_j$ of such alternatives are $\mathcal{S}_2=\mathcal{S}_3=\{yes, no\}$ and $\mathcal{S}_1=\{10, 20, 25, 40, 50\}\ \mu$g Hg/g of blood in the stool sample. Let $(g,k)$ with $g\in \mathcal{G}=\{F, M\}, k\in\{1, 2, 3, 4, 5\}$ be the target segment, where $F$ and $M$ refer to female and male, respectively, and numbers $k$ to the screening period. These periods correspond to a participant’s age so that period 1 refers to 60-year-olds, period 2 to 62-year-olds etc. The decisions about which alternatives $s_j$ (FIT cut-off level, incentive and invitation) to select for each target segment are modeled as binary decision variables $z_{g,k}(s_j) \in \{0,1\}$ so that $z_{g,k}(s_j)=1$ if and only if alternative $s_j$ is selected for segment $(g,k)$. A full screening strategy $Z$ is a collection of these decision variables for all target segments and can be written as $Z = \bigcup_{g,k} Z_{g,k}$ where the segment-specific strategies are defined as $Z_{g,k} = \bigcup_{j\in D} z_{g,k}(s_j), \forall g, k$.

Chance nodes correspond to returning a FIT sample, FIT results, continued participation, the discovery of polyps and polypectomy in a colonosopy, and adverse effects from the colonoscopy. In particular, nodes 5 and 8 reveal information about the participant's bowel state. Here, we assume that the bowel state $b \in \mathcal{B} = \{\mathrm{N, B, L, R}\}$  of a participant can be Normal (N), Benign growth (B), Large growth (L) or CRC (R). These bowel states are reflected at the population level by prevalences $\psi_{g,k,b}=N_{g,k,b} / N_{g,k}$, where $N_{g,k,b}$ is the number of participants in segment $(g,k)$ with  bowel state $b$, and $N_{g,k}$ is the total number of participants in segment $(g,k)$. The starting prevalences $\psi_{F,1,b}, \psi_{M,1,b}, \forall b\in \mathcal{B}$ as well as the segment-specific conditional probabilities for chance nodes are obtained from the literature (please see the Supplementary material for details). 

The progression of the prevalences of different bowel states are affected by two factors: natural progression and screening. The natural progression of colorectal cancer is reflected by transition probabilities $\mathbb{T}_{b,b'}^{g,k}$, i.e. the probability that the bowel state of the participant of sex $g$ is $b'$ in period $k+1$ given that at $k$ it was $b$. We assume that this progression follows the adenoma-carcinoma sequence, meaning the transition through bowel states can be represented by a linear recurrence relation. Screening, on the other hand, helps decrease the prevalences of abnormal states in the population, depending on the selected FIT cut-off level and incentive. In particular, we assume that any benign growth, large growth, or CRC found during the screening pathway is removed and that the bowel returns to a normal state. 
Taken together, the starting prevalences in screening period $k+1$ can be computed using the following difference equations: 

\begin{align}
    \psi_{g,k+1,B}(Z_{g,k}) &= (\psi_{g,k, B} - \tilde{\psi}_{g,k,B}(Z_{g,k}) )(1-\mathbb{T}^{g,k}_{\mathrm{B,L}}) +  \psi_{g,k,\mathrm{N}}\mathbb{T}^{g,k}_{\mathrm{N,B}} \label{eq:B_update}\\
    \psi_{g,k+1,\mathrm{L}}(Z_{g,k}) &= (\psi_{g,k, \mathrm{L}} - \tilde{\psi}_{g,k,\mathrm{L}}(Z_{g,k}) )(1-\mathbb{T}^{g,k}_{\mathrm{L,R}}) + (\psi_{g,k, \mathrm{B}} - \tilde{\psi}_{g,k,\mathrm{B}}(Z_{g,k}) )\mathbb{T}^{g,k}_{\mathrm{B,L}} \\
    \psi_{g,k+1,\mathrm{R}}(Z_{g,k}) &= \psi_{g,k, \mathrm{R}} - \tilde{\psi}_{g,k,\mathrm{R}}(Z_{g,k}) + (\psi_{g,k, \mathrm{L}} - \tilde{\psi}_{g,k,\mathrm{L}}(Z_{g,k}) )\mathbb{T}^{g,k}_{\mathrm{L,R}} \\
    \psi_{g,k+1,\mathrm{N}}(Z_{g,k}) &= 1 - \sum_{b \in \{\mathrm{B,L,R}\} } \psi_{g,k+1,b}(Z_{g,k})\quad , \label{eq:N_update}
\end{align}
where $\tilde{\psi}_{g,k,\mathrm{B}}(Z_{g,k})$ stands for the fraction of participants found to have bowel state $b$ in screening period $k$ as a result of applying screening strategy $Z_{g,k}$.

Following the notation in \citet{salo2022decision}, we define an \textit{information set} $\mathcal{I}(j)$ as the direct predecessors of a given node $j$, i.e. $\mathcal{I}(j) = \{i \in \mathcal{N} | (i,j) \in \mathcal{A}\}$. For instance, in Figure \ref{fig:combined_algo} the information set of Stage 5 (FIT result) consists of Stage 1 (FIT cut-off level), and Stage 4 (usable sample is returned). \textit{Information states} $s_{\mathcal{I}(j)} \in \mathcal{S}_{\mathcal{I}(j)}$ are defined as the combinations of states within all nodes of the the information set $\mathcal{I}(j)$. Relating to our previous example, the information states of Stage 5 are given as $s_{\mathcal{I}(5)} \in \mathcal{S}_{\mathcal{I}(5)} = \prod_{i \in \mathcal{I}(5)} \mathcal{S}_i = \mathcal{S}_1 \times \mathcal{S}_4$. For chance nodes, the conditional probability outcome $s_j$ depends on these information states. Specifically, the conditional probability of $s_j \in \mathcal{S}_j$ (where $j \in \mathcal{C}$) occurring can be given as $\mathbb{P}[X_j = s_j | X_{\mathcal{I}(j)} = s_{\mathcal{I}(j)}]$, where $X_{\mathcal{I}(j)}$ are the realizations of random variables $X_i$ in the information set $i \in \mathcal{I}(j)$. For example, in Figure \ref{fig:combined_algo} the conditional probabilities for the chance alternatives if returning a usable sample (Stage 4) are dependent on whether a person is invited to screen (Stage 3) and if they are offered an incentive (Stage 2). Please see \ref{AB:con_probs} for the conditional probability tables of each chance node in Figure \ref{fig:combined_algo}. Whereas, for decision nodes, the information states detail the information available when making the decision. For example, in Stage 7, when making the additional examination decision, the FIT result (Stage 5) and whether participation in the programme has continued (Stage 6) is available in the information state.

A \textit{local decision strategy} $Z_j$ for node $j \in \mathcal{D}$ is a function that maps each information state to a decision, i.e. $Z_j : \mathcal{S}_{\mathcal{I}(j)} \mapsto \mathcal{S}_j$. Binary decision variables $z(s_j|s_{\mathcal{I}(j)}) \in \{0,1\}$ denote the decision alternative $s_j \in j \in D$ that defines the local decision strategies. A \textit{global decision strategy} $Z_{g,k}$ (for a segment $g,k$) is a set of local decision strategies that contains one local strategy for each decision node. An example of a global decision strategy for the ID in Figure \ref{fig:combined_algo} could be $Z_{g,k} = \{Z_1 = \textit{50 ,} Z_2 = \textit{Yes ,} Z_3 = \textit{Yes}$.

Unlike a decision strategy that solely includes decisions, a \textit{path} $s \in \mathcal{S} = \bigtimes_{j=1}^N S_j$ is defined as the combination of all chance and decision node alternatives $s_j \in S_j \ \forall j \in C \cup D$, in which a global strategy and the value nodes are uniquely defined. For example, one possible path is: 
\begin{center}
  $s=$(75, No, Yes, Yes, Positive, Yes, Colonoscopy, Benign Adenoma, Polypectomy, None).
\end{center} 
The above path describes a sequence of events in which a person is invited ($s_3=$\emph{Yes}) to the screening program without monetary incentive ($s_2=$\emph{No}) using a FIT with a cut-off level value of 75 ng/ml ($s_1=75$). The invitee returns a usable sample ($s_4=$\emph{Yes}) that tests over 75 ng/ml and is thus scored as positive ($s_5=$\emph{Positive}), and contact is established with the local nurse ($s_6=$\emph{Yes}). A further examination is chosen to be a colonoscopy ($s_7=$\emph{Colonoscopy}), the result of which indicates a benign adenoma in the bowel ($s_8=$\emph{Benign Adenoma}). During the colonoscopy, the growth is removed from the bowel ($s_9$=\emph{Polypectomy}) and no adverse event occurs ($s_{10}$=\emph{None}). The scenario path defines the direct costs related to the participant's screening process, together with the detected cancer or growths, and whether a colonoscopy was administered.

The probability of a scenario path $\pi(s)$ occurring given the chances and the decisions made during the screening process are modeled as auxiliary variables. If the scenario path is not compatible with the decision strategy ($Z_{g,k}$) then $\pi(s) = 0$. For example, if on scenario path $s'$ the decision alternatives are $s_1 = 50, s_2 =$\textit{No}, but the local decision strategy has binary values that correspond to decision alternatives $s_1 = 75, s_2 =$\textit{No}, then $\pi(s')=0$, as the first decision for $s'$ does not match the current strategy. Otherwise, the path probability is equal to the upper limit, i.e.
\begin{equation}
    \pi(s) = \begin{cases}p(s) = \prod\limits_{j \in C}\mathbb{P}(X_j = s_j | X_{\mathcal{I}(j)}=s_{\mathcal{I}(j)}), \qquad \text{ if } Z_j(s_{\mathcal{I}(j)})=s_j \forall j\in \mathcal{D} \label{eq:p}\\
    0,\hspace{6.2cm} \text{ otherwise,}
    \end{cases}
\end{equation}
where $X_j$ is the random variable associated with chance node $j \in C$ and $p(s)$ denotes the upper bound for $\pi(s)$.


The optimization problem defined by an influence diagram is both discrete and linear. Thus, the set of all Pareto efficient (non-dominated) solutions can be calculated using an \emph{Modified Augmented Weighted Tchebychev} (MAWT) norm approach. The original algorithm with proofs is presented in \citet{Holzmann2018SolvingScalarizations}. We can denote the whole MAWT process as $(\Xi^*, \Phi^*) = \mathrm{MAWT}(ID(sex, k, {\psi}_{k}^{e}))$ where ${\psi}_k^{e}$ is the starting prevalence distribution. In the following model description, the solutions will all be dependent on the MAWT parameters, i.e. $\xi^* = \xi_{k, {\psi}^e_k}^{sex, *}$ and similarly for all components and sets. For the sake of avoiding notation clutter, we've excluded these indices from the discussion where possible without sacrificing clarity.

The approach relies on minimizing the MAWT norm generated by the algorithm. The \mbox{norm $||\cdot||^{\mathbf{w},\epsilon}$} selects the individual element with the highest value from the objective vector. 
The original MOMILP problem defined by the ID, utopian vector $\phi^{utopia}$, and nadir vector $\phi^{nadir}$ are given as inputs for the algorithm. The utopian vector is a weak lower bound for the Pareto optimal objective function values, which is determined by maximizing each objective individually. The nadir vector is a strong upper bound over the non-dominated solution set $\Phi$, assuming the objectives are minimised. The weights $\mathbf{w}$ and parameter $\epsilon$ are computed as in \citet{Holzmann2018SolvingScalarizations}.  

We can now formulate the optimisation problem solved by the MAWT algorithm. This  formulation captures the segment specific screening problem description above as a Multi-objective mixed integer linear problem (MOMILP). The variable $\mu$ (equation \eqref{eq:MAWT_obj}) is the MAWT norm, which is minimised on each iteration round and constraints \eqref{eq:MAWT_constr} are the actual objectives of the Decision Programming model. Equations \eqref{eq:z_sum}-\eqref{eq:pi_def} correspond to the Decision Programming formulation. 
\begin{align}
    &\min. \qquad \mu \label{eq:MAWT_obj}\\
    &\text{subject to } \\
     &\mu  \geq w_O|O - \phi^{utopia}_O| + \epsilon\sum_{O} w_O|O - \phi^{utopia}_O|, & \forall \  O\in \{\tilde{\Omega}, \tilde{N^{Col}}, \tilde{\psi}_B, \tilde{\psi}_L, \tilde{\psi}_R\} \label{eq:MAWT_constr}\\
     & \tilde{\Omega} = \sum_{s} \pi(s)U_{\mathrm{Cost}}(s)\\
     & \tilde{N}^{Col} = \sum_{s} \pi(s)U_{\mathrm{P Col.}}(s) \\
     & \tilde{\psi}_B = \sum_{s} \pi(s)U_{\mathrm{B}}(s) \\
     & \tilde{\psi}_L = \sum_{s} \pi(s)U_{\mathrm{L}}(s) \\
     & \tilde{\psi}_R = \sum_{s} \pi(s)U_{\mathrm{R}}(s) \\
     &\sum_{s} z(s_j | s_{\mathcal{I}(j)})  = 1, & \forall\ j \in \mathcal{D},\ s_{\mathcal{I}(j)} \in \mathcal{S}_{\mathcal{I}(j)} \label{eq:z_sum}\\
    &0  \leq \pi(s) \leq p(s), \qquad &\forall\ s\in \mathcal{S} \label{eq:pi_lim}\\
    &\pi(s)  \leq z(s_j | s_{\mathcal{I}(j)}), \qquad &\forall\ s\in \mathcal{S}, j \in \mathcal{D}\label{eq:pi_z}\\
   &\pi(s)  \geq p(s) + \sum_{j\in \mathcal{D}} z(s_j | s_{\mathcal{I}(j)}) - |\mathcal{D}|, &\forall\ s\in \mathcal{S}\label{eq:pi_lb}\\
    %
   %
    %
    %
    &z(s_j | s_{\mathcal{I}(j)}), \ U_i, \in \{0, 1\}, &\forall \ j \in D,\ s \in \mathcal{S}, i \in \{\mathrm{P Col., B, L, R }\} \label{eq:bin_domains}\\
    &\pi(s) \in \mathbb{R}.  &\forall\ s\in \mathcal{S} \label{eq:pi_def} 
\end{align}
Here $z(s_j | s_{\mathcal{I}(j)})$ are the actual decision variables that capture the segment specific strategy. The same strategy is reflected in the values of path probability variables $\pi(s)$. Values of $U_i(s)$ correspond to the utility node outcomes for the different paths in the ID. These are defined in \ref{s:utility_nodes}.  
When executed, the MAWT algorithm produces the complete efficient solution set $\Xi^* = \{\xi^* | \xi^* \succeq \xi,\ \xi^*,\xi \in \Xi \}$, where $\xi = (Z_{g,k}, \pi(s)) 
$ (i.e. a tuple of Pareto optimal variable values) and the non-dominated objective value set $\Phi^* = \{\tilde{\Omega}(\xi^*), \tilde{N}^{col}(\xi^*), \tilde{\psi}_{\mathrm{B}}(\xi^*),\tilde{\psi}_{\mathrm{L}}(\xi^*),\tilde{\psi}_{\mathrm{R}}(\xi^*) | \xi^* \in \Xi^*\}$



\subsection{Phase 1} 

In each period of Phase 1 the following steps are taken: (1) Sex- and period-specific age prevalences of CRC and benign and large adenomas are used as parameters in the optimization model. (2) The optimization model minimizes expected costs and expected colonoscopies while maximizing found cancers and growths in the bowel for that period. (3) The set of Pareto efficient screening strategies (i.e. paths through the influence diagram) for the specified period and sex are identified. (4) This set of efficient strategies is combined with the previous period's efficient set if it exists. (5) Any strategies that become inefficient as a result of the merging of sets are removed. (6) The sex- and period-specific prevalences are updated based on objective function values, that is, the number of found CRC and benign and large adenomas, and the natural transition of bowel states in the period interval. These updated prevalences are used in step (1) of the subsequent period. Steps (1) to (6) are repeated until all periods are optimized. 

\subsection*{Algorithm}

The algorithm used to solve the problem is presented in pseudocode in Algorithm \ref{alg:CRC_algo}.
\IncMargin{1em}
\SetAlCapSkip{2ex}
\begin{algorithm}
    \SetKwInOut{Input}{Inputs}\SetKwInOut{Output}{Output}
    \SetKwFunction{MAWT}{MAWT}\SetKwFunction{UpdateP}{UpdatePrevalences}\SetKwFunction{TotalP}{TotalPrevalences}\SetKwFunction{RemoveDom}{RemoveDominated}
    \SetKwBlock{Initialize}{Initialize}{end}
    \Input{Process model as Influence Diagram (ID), $\psi_{F,1,b}, \psi_{M,1,b}, \forall b\in\mathbf{B}$, screening periods $k\in \{1\dots K\}$ } 
    \Output{The optimal strategy  $Z^* \in \mathcal{E}_{F,K} \times \mathcal{E}_{M,K}, \mathcal{E}_{g,k} = \{e| e \succeq e' \}$ for each sex $g$ , and the related objective value $\Psi_{\mathrm{R}}(Z^*)$.}
    $\mathcal{E}_{F,k}, \mathcal{E}_{M,k} = \varnothing \quad \forall k \in [1 \dots K]$ \tcp*[r]{Initialize set of Efficient solutions as empty sets.}

    \For{$g \in \{\mathrm{F, M}\}$}{
        \label{alg1:sexloop}\tcp{First period, $k=1$, handled separately.}
        
        \tcp{Compute Pareto efficient solutions for first period}
        
        $\mathcal{E}_{g,1} \leftarrow$  \MAWT $(ID(g,1), \bar{\psi}_{g,1}))$ \tcp*[r]{$\bar{\psi}_{g,k} = [\psi_{g,k,b}]_{b \in \mathcal{B}}$} \label{alg1:first_E}

        $\bar{\psi}_2^{e^{g}_1} \leftarrow \UpdateP(e^{g}_1, \bar{\psi}_{g,1}) \quad \forall e_{g,1} \in \mathcal{E}_{g,1}$ \label{alg1:first_update}\;
        
        \label{alg1:first_col}$N^{col}_{total} \leftarrow N^{col}(e_{g,1})$

        \For{$k \in \{2,\dots,K\}$}{\label{alg1:periodloop}
            \For{$e \in \mathcal{E}_{g,k-1} $}{
            \tcp{Compute Pareto efficient solutions for period $k$ based on preceding strategy $e$}$\mathcal{Z}_{g,k}^*  \leftarrow$ \MAWT $(ID(g,k, \bar{\psi}(e)))$\; 
            
            \For{$Z_{g,k} \in \mathcal{Z}_{g,k}^*$}{
            
                $e_k \leftarrow (e,Z_{g,k})$ \tcp*[r]{Extend $e$ by $Z_{g,k}$ into a new strategy $e_k$, covering steps until $k$}
                $\mathcal{E}_{g,k} \leftarrow \mathcal{E}_{g,k} \bigcup e_{g,k} $\label{alg1:z_ext} \tcp*[r]{Collect found strategies for step k.}
                
                
               $\bar{\psi}_{g,k}^{e_{g,k}} \leftarrow \UpdateP(e_{g,k}, \bar{\psi}_{g,k-1})$ \label{alg1:k_update}\;
                
                $\bar{\Psi}_k^{e_{g,k}} \leftarrow $ \TotalP($\bar{\psi}_{g,k}^{e_{g,k}}, \bar{\Psi}_{k-1}^{e}$)\;
                
                $N^{col,e_{g,k}}_{total} \leftarrow N^{col,e}_{total}  + \tilde{N}^{col}(Z_{g,k})$ \label{alg1:colon_sum}\tcp*[r]{Add new colonoscopies from $Z_{g,k}$ to preceding colonoscopies.}
                
                \If{$N^{col,e_{g,k}}_{total}\geq N^{col,max}$\label{alg1:colon_constr}} 
                    {
                    $\mathcal{E}_{g,k} \leftarrow \mathcal{E}_{g,k} \cap e_{g,k}$
                    }
                }
            }
            $\mathcal{E}_{g,k} \leftarrow \RemoveDom(\mathcal{E}_{g,k})$ \label{alg1:rem_dom}
        }
    }
    $\mathcal{E}_K \leftarrow \RemoveDom(\mathcal{E}_{F,K} \times \mathcal{E}_{M,K})$
    \caption{Phase 1 algorithm}\label{alg:CRC_algo}
\end{algorithm}\DecMargin{1em}
As main inputs, the algorithm requires 1) Influence Diagram (ID) describing the screening process per one population segment defined by age and sex 2) MAWT algorithm for computing a Pareto front for one ID 3) Starting prevalences for the first age segments (one prevalence vector per sex) 3) a prevalence update function that takes as input a prevalence at period $k$ and a period specific screening strategy, and outputs an updated prevalence vector for $k+1$ 5) a function for combining age group specific prevalences into a total population prevalence. The output is a Pareto efficient set of sex specific screening strategies, i.e. combinations of period-specific strategies and their related objective values and supportive information.

The process starts by initialising empty sets of efficient solutions (one per sex). The sexes are handled separately (line \ref{alg1:sexloop}). The algorithm then proceeds to compute strategies, one for each screening segment. The first period is handled separately, as the process is simpler than for the remaining periods. A set of Pareto efficient, period specific solutions is computed by iteratively solving the MAWT version of the Influence Diagram (see section \ref{sec:ID_MAWT}) until the complete Pareto front is found. These strategies are saved into efficient solution set $\mathcal{E}_{g,1}$ on line \ref{alg1:first_E}. One element in this set contains a strategy (only for the first period at this point) and the related objective function values. On line \ref{alg1:first_update} expected prevalences for bowel states after screening and 2 years of waiting for the next period, $\bar{\psi}^{z_1}_1$, are computed for each found strategy. These will act as starting points for period 2 optimization. 

Next, a similar process is performed for periods $k\in [2,\dots,K]$ using the found solutions $e \in \mathcal{E}^{sex}_{k-1} $ as inputs. For each surviving previous period solution, a new Pareto front is created by solving the MAWT problem, now with a prevalence $\bar{\psi}^e_{k}$  based on the previous period efficient strategy $e$ as a starting point. Each of the found strategies is added to the strategy included in $e$ (row \ref{alg1:z_ext}), creating a history of decisions leading up to this decision. Again, the prevalences are updated as on the first period. An additional step now computes the total population prevalence for the population consisting of sex specific age groups corresponding to periods $[1, \dots, k]$. The cancer component in this total prevalence vector is used as one of the objectives in filtering out Pareto non-efficient strategies, together with cancer and large adenoma prevalences of the current age group.
If the sum of colonoscopies (computed on row \ref{alg1:colon_sum}) at this period, with found strategy, does not exceed the maximum number of coloscopies (row \ref{alg1:colon_constr}), the solution is added to the set of found solutions.

Finally, on row \ref{alg1:rem_dom} the solutions found at period $k$ are compared in terms of total cancer prevalence, period $k$
age group cancer and large adenoma prevalences and required total colonoscopies, and the Pareto inefficient solutions are removed. Including the different prevalences at this point ensures that on the next round, the starting points are efficient in terms of the overall objective of minimising cancer prevalence in the total population but also that the starting point on the next round $k+1$ is efficient. Large adenomas are included also in the starting point comparison, as the influence the generation of new cancers at $k+1$. The algorithm finishes when the predetermined period $K$ has been reached.

\subsection{Phase 2}
In this phase, a strategy for each age of males and females is selected in order to minimize CRC prevalence in the entire screening population with respect to a colonoscopy constraint. In each period of Phase 1 the efficient solutions of the optimization problem \eqref{eq:MAWT_obj} - \eqref{eq:pi_def} are combined with the previous period's solution set with any now inefficient solutions removed. This means that in Phase 2, only the solution set from the final period is relevant as it contains the optimal solutions from all periods. Without a temporal aspect, the Phase 2 optimization problem is simple and straightforward to solve.

For this problem we take a utilitarian approach, in that we wish to minimize the population-level cancer prevalence subject to a population-level resource constraint. To achieve this, we solve the binary linear programming model: 
\begin{align}
     \min_{\bar{x}_{\mathrm{F},j}, \ \bar{x}_{\mathrm{M},j}}. \quad & \Psi_R = \frac{N^R_{\mathrm{F},j} \bar{x}_{\mathrm{F},j} + N^R_{\mathrm{M},j} \bar{x}_{\mathrm{M},j}}{N_\mathrm{F} + N_\mathrm{M}},\\
     \text{subject to} \quad & N_\mathrm{F} N^{col}_{\mathrm{F},j} \bar{x}_{\mathrm{F},j} + N_\mathrm{M} N^{col}_{\mathrm{M},j} \bar{x}_{\mathrm{M},j} \leq N^{Col,Max},\\
     & \bar{x}_{\mathrm{M},j} \in \{0,1\}, \qquad \forall j \in \{1,\dots,J_M\}\\
     & \bar{x}_{\mathrm{F},j} \in \{0,1\}, \qquad \forall j \in \{1,\dots,J_F\}\\
     & \sum_j x_{g,j} = 1 \qquad \forall g \in \{\mathrm{F}, \mathrm{M}\}\,
\end{align}
where $N$ refers to the population size; $\mathrm{M}$ and $\mathrm{F}$ correspond to males and females; $i$ is the bowel state colorectal cancer; $j$ is a strategy from the set of efficient strategies, sized $J_M$ and $J_F$ for males and females respectively, identified in the final period of Phase 1; $N^{col}_{\mathrm{M},j}$ and $N^{col}_{\mathrm{F},j}$ are the expected number of colonoscopies in strategy $j$ for males and females, respectively; and, $N^{col,max}$ is the population-level colonoscopy budget constraint. The binary variables $x_{g,j}$ are used to select one strategy for both sex.

\section{Results and policy recommendations} \label{Sec:Results}

In this section, we present optimal screening policies for two situations (i) when decisions regarding monetary incentives for participation are not available, and (ii) when these incentives decisions are available. Firstly, the situation in which no incentives are available reflects the current Finnish CRC screening programme. These results will provide optimal screening policy recommendations to aid decision makers when considering the choice of FIT cut-off level for different age and sex groups. Additionally, the expected costs of each policy will be presented. To the best of our knowledge, no other cost analysis for the current Finnish CRC screening programme has been published in the public domain. Secondly, results and recommendations for a hypothetical situation in which monetary incentives to participate in screening are presented, along with their expected costs.

To assess the policy and clinical implications of different colonoscopy resource constraints are compared to the baseline case of no screening. To also study the role of the available colonoscopy resources, we solve the problem in 4 cases, i.e. for 8000, 12000, 16000 and 20000 colonoscopies as the maximum colonoscopy limit for all 5 groups screened.





\subsection{Results for no monetary incentive option}

\begin{table}[!ht]
\centering
\begin{tabular}{|l|l|ccccc|}
\hline
Case (max colonoscopies) & Sex & Age 60 & Age 62 & Age 64 & Age 66 & Age 68 \\ \hline
Case 1 (8 000) & Male & 50 & 50 & 50 &  50 &  50 \\ 
 & Female & - & 40 &  25 & 25 & 25 \\ \hline
Case 2 (12 000)& Male & 50 & 25 & 50 &  20 &  20 \\ 
 & Female & 20 & 25 & 20 & 20 &  10 \\ \hline
Case 3 (16 000)& Male & 20 & 20 &  20 & 20 & 20 \\ \
 & Female &  10 &  10 &  10 &  10 &  10 \\ \hline
Case 4 (20 000)& Male &  10 &  10 &  10 &  10 &  10 \\ 
 & Female &  10 &  10 &  10 &  10 &  10 \\ \hline
\end{tabular}
\caption{Optimized screening policies with no option of incentive: FIT cut-off level in $\mu$g/g, no invite represented by -, colonoscopy performed when positive FIT.}
\label{tab:policy_no_inc}
\end{table}

In each case and screening age, the FIT cut-off level used is presented in Tables 1 and 2. For example, for Case 2 with no option of incentive (Table 1), the male cut-off level is 50$\mu$g/g and the female is 20$\mu$g/g for age 60. Whereas, the cut-off level of 250$\mu$g/g is used for both males and females aged 62. 

\begin{figure}[!ht]
    \centering
    \includegraphics[width=0.9\textwidth]{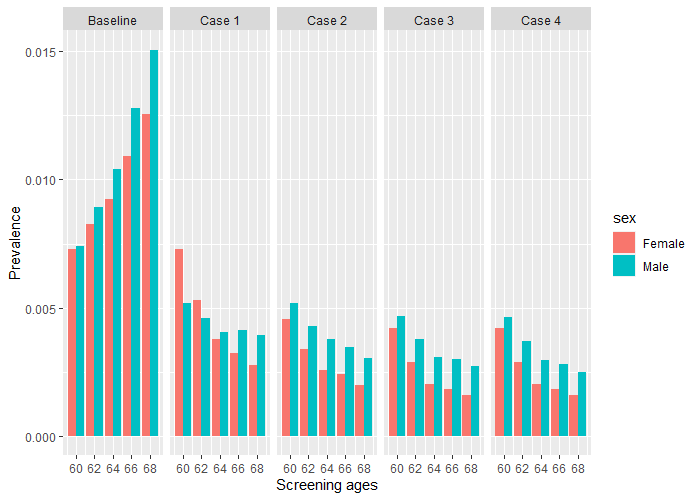}
    \caption{Cancer prevalence over screening rounds with no option of incentive, policies given in Table \ref{tab:policy_inc}}
    \label{fig:CRCprev_noCVaR_no_inc}
\end{figure}

\subsection{Results for monetary incentive option}

\begin{figure}[ht!]
    \centering
    \includegraphics[width=0.9\textwidth]{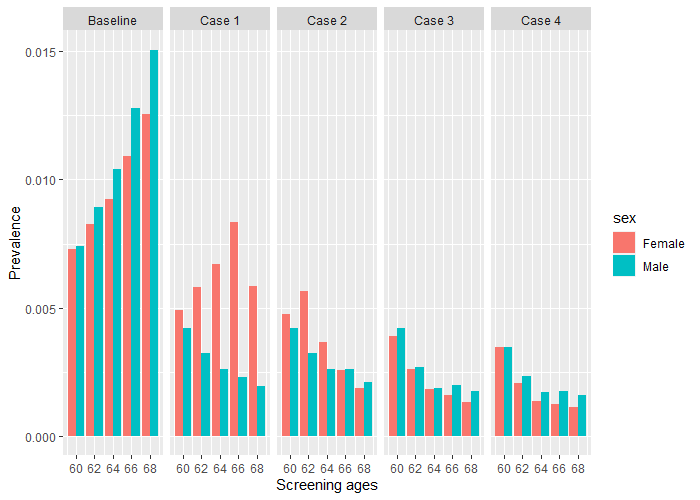}
    \caption{Cancer prevalence over screening rounds with option of incentive, policies given in Table \ref{tab:policy_inc}}
    \label{fig:CRCprev_noCVaR_inc}
\end{figure}

The optimal strategies for different cases are in Table \ref{tab:policy_inc} and the corresponding cancer prevalences in Figure \ref{fig:CRCprev_noCVaR_inc}. In our model the use of incentives can be beneficial, even at the cost of lower sensitivity due to higher cut-off levels, which can be seen by comparing these strategies to the ones where incentives were not available. The benefits of incentives are apparent in Figure \ref{fig:CRCprev_noCVaR_inc} where the total prevalences are lower for the incentivized strategies, except for case 1 where the availability of incentives caused the screening to focus heavily on males. This implies that ways of increasing the participation rate of males could prove valuable for total cancer prevalence in the population, regardless of the exact method for achieving such increase.
\begin{table}[ht!]
\centering
\begin{tabular}{|l|l|ccccc|}
\hline
Case (max colonoscopies) & Sex & Age 60 & Age 62 & Age 64 & Age 66 & Age 68 \\ \hline
Case 1 (8 000)  & Male & 50+i & 50+i & 50+i &  20+i &  20+i \\ 
 & Female & 50 & - &  - & - & 25 \\ \hline
Case 2 (12 000)& Male & 50+i & 50+i &  50+i &  40+i &  20+i \\ 
 & Female & 40 &  - &  40+i & 10+i &  10+i \\ \hline
Case 3 (16 000)& Male & 50+i & 20+i & 20+i &  25+i &  20+i \\ 
 & Female &  20+i &  25+i &  25+i &  10+i &  10+i \\ \hline
Case 4 (20 000) & Male & 20+i & 20+i &  20+i &  20+i &  10+i \\ 
 & Female &  10+i &  10+i &  10+i &  10+i &  10+i \\ \hline
\end{tabular}
\caption{Optimized screening policies with option of incentive (+i): FIT cut-off level in $\mu$g/g, no invite represented by -, colonoscopy performed when positive FIT.}
\label{tab:policy_inc}
\end{table}

\section{Discussion} \label{Sec:Diss}


In this manuscript we have proposed a novel approach for optimizing a cancer screening program and applied it to the Finnish colorectal cancer screening program. The optimization approach combines several different algorithms and multi-objective optimization to find an optimal solution to a multi-period screening process captured in IDs, and respective of a colonoscopy constraint. Given the ID and its related parameter values and a sex-specific starting prevalences, our approach is able to find an optimal strategy for a 5-period screening process from among approximately 26 billion possible alternative strategies. Compared to previous approaches to screening program optimization documented in literature, one of the main strengths of our approach is the ability to optimize the strategy over several screening rounds.

We computed results for 4 different colonoscopy resource levels. The different result sets show a clear differences in selected FIT test cut-off levels and prioritization between male and female participants, and the potential benefits from varying the FIT cut-off level between different sex and age segments, when resources are limited. We also present separate results for cases where a hypothetical incentive is available, or not. There is a clear difference in the incentivized results compared to non-incentivized in the low-resource cases. With incentives available, the results emphasize testing males with incentives at the cost of screening females. 

This study opens up several avenues for future research. Modeling-wise, we have made some simplifications, on parts of which the model could be improved. In particular our assumptions regarding the transition probabilities and healing of participants for whom an abnormal bowel was discovered are relatively strong. Improvement of the transition model could therefore lead to more accurate estimates and strategy proposals. Inclusion of mortality as an objective could also provide important insights into the outcomes of the model. Due to lack of data, mortality was outside the scope of this study but the situation is likely to improve as the screening program runs for more years. The approach we have presented is flexible in the sense that it is also capable of solving the problem under different types of constraints. 

Our approach is computationally quite intensive, as it, in practice, generates a multi-periodic decision tree to be solved. However, due to the nature of the algorithm, the main memory requirement is set by the structure of the influence diagram, i.e. the period-specific model, whereas the addition of periods increases only the computation time. Separating the branches of the whole multi-period tree into separate sub-problems opens up the possibility to solve the in parallel, which could be used to decrease the required computation time. Additional problem-specific constraints can possibly also be used to decrease computational time.


\section{Acknowledgements}

The calculations presented above were performed using computer resources within the Aalto University School of Science “Science-IT” project.

This research has been supported by The Foundation for Economic Education (grant number 16-9442) and The Paulo Foundation. 

\bibliographystyle{elsarticle-harv}
\biboptions{authoryear}

\bibliography{references} 




%
%
%
\newpage
\appendix

\section{Conditional probabilities and values in the influence diagram}\label{AB:con_probs}

\subsection{Conditional probabilities for stage 4: Usable sample returned}
The conditional probabilities for returning a sample in good condition are
\begin{equation}
    \mathbb{P}(X_4 = s_4 | X2=s_2, X_3=s_3), \quad  \forall s_4 \in \{\mathrm{No, Yes}\}, s_2 \in \{\mathrm{No, Yes}\}, s_3 \in \{\mathrm{No, Yes}\} \label{eq:C4}.
\end{equation}
This can be presented as a matrix:

$\mathbb{P}(X_4 = s_4 | X2=s_2, X_3=s_3) = $
  \begin{blockarray}{*{2}{c} cc}
    \begin{block}{*{2}{>{\footnotesize}c<{} >{\footnotesize}c<{}} }
      No & Yes & $s_2$ & $s_3$ \\
    \end{block}
    \begin{block}{[*{2}{c}]>{\footnotesize}c<{} >{\footnotesize}c<{} }
      1.0 & 0.0 & No & No\\
    $1- P^{ret. OK}_{g,a}$ & $P^{ret. OK}_{g,a}$ & No & Yes \\
    1.0 & 0.0 & Yes & No\\
    $(1- P^{ret. OK}_{g,a})/2$ & $P^{ret. OK}_{g,a} + (1- P^{ret. OK}_{g,a})/2$ & Yes & Yes \\
    \end{block}
  \end{blockarray}

We have assumed that incentivizing the invited persons has the effect of halving the probability of not returning a sample.
The probability $P^{ret. OK}_{g,a}$ is calculated as
\begin{equation}
    P^{ret. OK}_{g,a} = P(\mathrm{return} | Age = a, Sex = g) P(\mathrm{Sample OK}).
\end{equation}

\subsection{Conditional probabilities for stage 5: FIT test result}

The conditional probabilities for chance node $C_5$ are 
\begin{equation}
    \mathbb{P}(X_5 = s_5 | X_1=s_1, X_4=s_4), \quad  \forall s_5 \in \{NA, -, +\}, s_1 \in \{50, 75, 100, 150, 250\}, s_4 \in \{No, Yes\} \label{eq:C5}
\end{equation}

This can be presented as a matrix:

  $\mathbb{P}(X_5 = s_5 | X_1=s_1, X_4=s_4) =$
  \begin{blockarray}{*{3}{c} cc}
    \begin{block}{*{5}{>{\footnotesize}c<{} } }
      NA & + & - & $s_1$ & $s_4$ \\
    \end{block}
    \begin{block}{[*{3}{c}]>{\footnotesize}c<{} >{\footnotesize}c<{} }
      1.0 & 0.0 & 0.0 &  50 & No\\
    0.0 & $P^{FIT+}_{50}$ & $1-P^{FIT+}_{50}$ & 50 & Yes \\
    1.0 & 0.0 & 0.0 & 75 & No\\
    0.0 & $P^{FIT+}_{75}$ & $1-P^{FIT+}_{75}$ & 75 & Yes\\
    1.0 & 0.0 & 0.0 & 100 & No\\
    0.0 & $P^{FIT+}_{100}$ & $1-P^{FIT+}_{100}$ & 100 & Yes \\
     1.0 & 0.0 & 0.0& 150 & No \\
    0.0 & $P^{FIT+}_{150}$ & $1-P^{FIT+}_{150}$ & 150 & Yes \\
     1.0 & 0.0 & 0.0& 250 & No \\
    0.0 & $P^{FIT+}_{250}$ & $1-P^{FIT+}_{250}$ & 250 & Yes \\
    \end{block}
  \end{blockarray}

where the positive result probabilities for FIT threshold levels $l$ $P^{FIT+}_{l}$ are calculated as follows:
\begin{align}
  &  P^{FIT+}_{l} = \sum_{b}P(FIT+ | l, b), \quad b \in \{\mathrm{Normal, Benign, Large, CRC}\} \\
  &  P(FIT+ | l, b) =  \sigma_{l, b}^{+, FIT} \psi_b, \quad b \in \{\mathrm{Benign, Large, CRC}\} \\
  &  P(FIT+ | l, \mathrm{Normal}) = (1-\sigma_{l, \mathrm{Normal}}^{-, FIT}) \psi_{\mathrm{Normal}}.
\end{align}
Here $\sigma^{+, FIT}_{l,b}$ is the sensitivity of the FIT test for detecting bowel states $b$ given threshold $l$ and $\sigma^{-, FIT}_{l,\mathrm{Normal}}$ is the corresponding specificity. The prevalences are segment (sex, age), period and preceding strategy specific, but for sake of avoiding index clutter, the information is omitted from the formulation.
For situations where no usable sample was returned, i.e. $s_4=\mathrm{No}$, the only possible outcome for the FIT test result is NA. 

\subsection{Conditional probabilities for stage 6: Continued participation}

The conditional probabilities for chance node $C_6$ are 
\begin{equation}
\mathbb{P}(X_6 = s_6 | X_5 = s_5), \quad  \forall s_6 \in \{No, Yes\}, s_5 \in \{\mathrm{NA}, +, -\} \label{eq:C6}
\end{equation}
which can be presented as a matrix

 $\mathbb{P}(X_6 = s_6 | X_5 = s_5) =$
 \begin{blockarray}{*{2}{c} c}
   \begin{block}{*{3}{>{\footnotesize}c<{} } }
     No & Yes & $s_5$  \\
   \end{block}
   \begin{block}{[*{2}{c}]>{\footnotesize}c<{} }
      1 & 0 & NA\\
      $1-P(Yes | Sex, Age)$ & $P(Yes | Sex, Age)$ & + \\
      1 & 0 & -\\
    \end{block}
 \end{blockarray}

\subsection{Conditional probabilities for stage 8: Examination result}

The conditional probabilities for chance node $C_8$ are 
\begin{equation}
\begin{array}{l}
    \mathbb{P}(X_8 = s_8 | X_1=s_1, X_6=s_6, X_7 = s_7), \quad \\ \forall s_8 \in \{NA, Normal, Benign, Large, CRC\}, s_1 \in \{50, 75, 100, 150, 250\},\\ s_6 \in \{No, Yes\}, s_7 \in \{No, Yes\} \label{eq:C8}
    \end{array}
\end{equation}

\newpage
This can be presented as a matrix:
  $\mathbb{P}(X_8 = s_8 | X_1=s_1, X_6=s_6, X_7 = s_7) = $\\
  \begin{blockarray}{*{5}{c} *{3}{c}}
    \begin{block}{*{8}{>{\footnotesize}c<{} } }
      NA & Normal & Benign & Large & CRC & $s_1$ & $s_6$ & $s_7$ \\
    \end{block}
    \begin{block}{[*{5}{c}]>{\footnotesize}c<{} >{\footnotesize}c<{} >{\footnotesize}c<{} }
      1.0 & 0.0 & 0.0 & 0.0 & 0.0 &  50 & No & No\\
      1.0 & 0.0 & 0.0 & 0.0 & 0.0 &  50 & No & Yes\\
      1.0 & 0.0 & 0.0 & 0.0 & 0.0 &  50 & Yes & No\\
      0.0 & $P^\mathrm{COL, Normal}_{50}$ & $P^\mathrm{COL, Benign}_{50}$ & $P^\mathrm{COL, Large}_{50}$ & $P^\mathrm{COL, CRC}_{50}$ & 50 & Yes & Yes\\
      1.0 & 0.0 & 0.0 & 0.0 & 0.0 &  50 & No & No\\
      1.0 & 0.0 & 0.0 & 0.0 & 0.0 &  50 & No & Yes\\
      1.0 & 0.0 & 0.0 & 0.0 & 0.0 &  50 & Yes & No\\
      0.0 & $P^\mathrm{COL, Normal}_{75}$ & $P^\mathrm{COL, Benign}_{75}$ & $P^\mathrm{COL, Large}_{75}$ & $P^\mathrm{COL, CRC}_{75}$ & 75 & Yes & Yes\\
      1.0 & 0.0 & 0.0 & 0.0 & 0.0 &  100 & No & No\\
      1.0 & 0.0 & 0.0 & 0.0 & 0.0 &  100 & No & Yes\\
      1.0 & 0.0 & 0.0 & 0.0 & 0.0 &  100 & Yes & No\\
      0.0 & $P^\mathrm{COL, Normal}_{100}$ & $P^\mathrm{COL, Benign}_{100}$ & $P^\mathrm{COL, Large}_{100}$ & $P^\mathrm{COL, CRC}_{100}$ & 100 & Yes & Yes\\
      1.0 & 0.0 & 0.0 & 0.0 & 0.0 &  150 & No & No\\
      1.0 & 0.0 & 0.0 & 0.0 & 0.0 &  150 & No & Yes\\
      1.0 & 0.0 & 0.0 & 0.0 & 0.0 &  150 & Yes & No\\
      0.0 & $P^\mathrm{COL, Normal}_{150}$ & $P^\mathrm{COL, Benign}_{150}$ & $P^\mathrm{COL, Large}_{150}$ & $P^\mathrm{COL, CRC}_{150}$ & 150 & Yes & Yes\\
      1.0 & 0.0 & 0.0 & 0.0 & 0.0 &  250 & No & No\\
      1.0 & 0.0 & 0.0 & 0.0 & 0.0 &  250 & No & Yes\\
      1.0 & 0.0 & 0.0 & 0.0 & 0.0 &  250 & Yes & No\\
      0.0 & $P^\mathrm{COL, Normal}_{250}$ & $P^\mathrm{COL, Benign}_{250}$ & $P^\mathrm{COL, Large}_{250}$ & $P^\mathrm{COL, CRC}_{250}$ & 250 & Yes & Yes\\
    \end{block}
  \end{blockarray}

where the probabilities $P^{COL, r}_{l}$ for colonoscopy results given a certain FIT level $l$  are calculated as follows:
\begin{align}
    & P^{\mathrm{COL}, r}_{l} = \sigma^\mathrm{+, COL}_{b} P(b | \mathrm{FIT+}) \\
    & P(b | \mathrm{FIT+}, l) = \frac{\sigma^\mathrm{+, FIT}_{l,b} \psi_b}{P(\mathrm{FIT+})} \\
    & P(\mathrm{FIT+}, l) = \sum_{b\ne \mathrm{Normal}} \sigma^\mathrm{+, FIT}_{l,b} \psi_b + (1-\sigma^\mathrm{-, FIT}_{l,\mathrm{Normal}}) \psi_\mathrm{Normal}.
\end{align}
Here we have assumed that the specificity of colonoscopies is perfect, i.e. $P( + | Normal ) = 0$. The sensitivities and specificities are acquired from data as are the prevalences of bowel states in the first period. In further periods, the prevalences are acquired from optimization results of preceding periods. The prevalences are segment (sex, age), period and preceding strategy specific, but for sake of avoiding index clutter, the information is omitted from the formulation.

\subsection{Conditional probabilities for stage 9: Polyp found}

The conditional probabilities for finding a polyp, and thus leading to a polypectomy, are
\begin{equation}
    \mathbb{P}(X_9 = s_9 | X8=s_8), \quad  \forall s_9 \in \{\mathrm{No\ result, No\  polyp, Polyp}\}, s_8 \in \{\mathrm{NA, Normal, Benign, Large, CRC}\} \label{eq:C9}.
\end{equation}
This can be presented as a matrix:

  $\mathbb{P}(X_9 = s_9 | X_8=s_8) =$
  \begin{blockarray}{*{4}{c}}
    \begin{block}{*{4}{>{\footnotesize}c<{}  } }
      {No\ result} & {No polyp} & {Polyp} & $s_8$ \\
    \end{block}
    \begin{block}{[*{3}{c}] >{\footnotesize}c<{} }
          1 & 0 & 0 & {NA}\\
          0 & 1 & 0 & {Normal}\\
          0 & 0 & 1 & {Benign}\\
          0 & 0 & 1 & {Large}\\
          0 & 0 & 1 & {CRC}\\
    \end{block}
  \end{blockarray}

This means that if any growth is found, the finding of a polyp is also assumed.

\subsection{Conditional probabilities for stage 10: Adverse event}

The conditional probabilities for an adverse event (bleed or a perforation) during polypectomy are 
\begin{equation}
    \mathbb{C}_{10} = \mathbb{P}(X_{10} = s_{10} | X9=s_9), \quad  \forall s_{10} \in \{\mathrm{None, Bleed, Perforation}\}, s_9 \in \{\mathrm{NA, No polyp, Polyp}\} \label{eq:C10}.
\end{equation}
This can be presented as a matrix:

  $\mathbb{C}_{10} =$
  \begin{blockarray}{*{4}{c}}
    \begin{block}{*{4}{>{\footnotesize}c<{}} }
      {No adverse event} & {Bleed} & {Perforation} & $s_9$ \\
    \end{block}
    \begin{block}{[*{3}{c}] >{\footnotesize}c<{} }
          1 & 0 & 0 & {No result}\\
          $1-P^\mathrm{Bleed}-P^\mathrm{Perf. (wo Poly.)}$ & $P^\mathrm{Bleed}$ & $P^\mathrm{Perf. (wo Poly.)}$ & No polyp\\
          $1-P^\mathrm{Bleed}-P^\mathrm{Perf. (w Poly.)}$ & $P^\mathrm{Bleed}$ & $P^\mathrm{Perf. (w Poly.)}$ & {Polyp}\\
    \end{block}
  \end{blockarray}

These probabilities are acquired from data.

\subsection{Conditional values for utility nodes}
\label{s:utility_nodes}
The conditional values for all growth types (Benign, Large, CRC), ${U}_{\mathrm{B}}, {U}_{\mathrm{L}}, {U}_{\mathrm{R}}$ are simply 1 if such a growth is found or 0 if not. Similarly, the value for colonoscopies, ${U}_{\mathrm{P Col.}}$, is -1 if a colonoscopy is performed and 0 otherwise. The negativity of the value is just a stylistic choice to reflect the 'resource consuming' role of performing a colonoscopy.

The cost value ${U}_{\mathrm{Cost}}$ is conditional to all cost generating stages. It can be formulated as
\begin{equation}
    \mathbb{U}_{14} = \Omega | (S_j = s_j) = \sum_{j} \Omega(s_j),\quad j \in \{2,3,4,7,8,9,10\}
\end{equation}
where the individual cost values for stage realizations are derived from data.

\end{document}